# USING FIBONACCI NUMBER TO INTEGRATE 2 × 2 AND 3 × 3 MATRICES

Honer N. Abdullah[*], Delbrin H. Ahmed and Muwafaq M. Salih
Dept. of Mathematics, College of Basic Education, University of Duhok, Kurdistan Region-Iraq

**ABSTRACT**
The aim of this paper is to make the clarification of images faster by the formula that Franciszekn made for matrices integrations and this made Sukhvinder's Algarithm complicate and slower. Tis paper uses Fibonacci number to determine integration formulas for matrices of order 2 and 3 in order to make the process of images clarification shorter.

## 1. INTRODUCTION

The Fibonacci numbers (Fibonacci sequence) were invented by Italian Leonardo Pisano Bigollo in his book called''Liber Abaci''. Furthermore, the integration of matrices is used to avoid image ambiguous. Franciszekn[1] used anti-orthogonality of matrix to find its integration, but that formula made the Sukhvinder's[7] algorithm complicate and the process became slow where the orthogonality founded interest in the [4],[7],[6]. Penner[3] use double integral to find the matrix integration. While Rnadal [4] found the derivative of the matrix for some special cases. In this paper, we used Fibonacci number formula to find the formula of integration of matrix in diminution 2 and 3.

## 2. SOME DEFINITIONS AND MATRIX DERIVATIVE

In this section, we see the basic definition that are going to be used in section two moreover some proposition of how to derivative matrix.

**Definition of Matrix 2.1[4]** The order rectangular $n = 1, 2, \ldots, i,\ m = 1, 2, \ldots, j$ where $b_{ij} \in R$ then

$$B = \begin{bmatrix} b_{11} & b_{12} & \ldots & b_{1m} \\ b_{21} & b_{22} & \ldots & b_{2m} \\ \vdots & \vdots & & \vdots \\ b_{n1} & b_{n2} & \ldots & b_{nm} \end{bmatrix}$$

is said to be matrix of $n \times m$ dimension.

**Definition of First Order Partial Derivative 2.2[4]**
Let $y = \alpha(x)$ where $y$ is vector with $n$ element and $x$ is $m$ element vector so

$$\frac{\partial y}{\partial x} = \begin{bmatrix} \frac{\partial y_1}{\partial x_1} & \ldots & \frac{\partial y_1}{\partial x_n} \\ \vdots & \ddots & \vdots \\ \frac{\partial y_m}{\partial x_1} & \ldots & \frac{\partial y_m}{\partial x_n} \end{bmatrix}$$

The matrix we get is $n \times m$ of first order partial derivative.

**Definition 2.3[3]** The Fibonacci number $F(k)$ are defined by the equation $F(k) = U_{k-1} + U_{k-2}$ where $U_1 = U_2 = 1$ and $k \geq 3$.

[*] E-mail: honer.naif@uod.ac

## 3. PROPOSITIONS

The following propositions as in[5]

**Proposition 3.1** Let $y = Ax$ where y is $m \times 1$, $x$ is $n \times 1$, $A$ is $m \times n$, and $A$ independent of $x$, then
$$\frac{\partial y}{\partial x} = A$$

**Proposition 3.2** Let $y = Ax$ where $y$ is $m \times 1$, $x$ is $n \times 1$, $A$ is $m \times n$, and $A$ independent of $x$, suppose that $x$ is a function of the vector $z$, while $A$ is independent of $z$. Then
$$\frac{\partial y}{\partial z} = A \frac{\partial x}{\partial z}.$$

**Proposition 3.3** Let the scalar α be defined by $\alpha = y^T A x$, where y is $m \times 1$, $x$ is $n \times 1$, $A$ is $m \times n$, and $A$ is independent of $x$ and $y$, then
$$\frac{\partial \alpha}{\partial x} = y^T A$$

And
$$\frac{\partial \alpha}{\partial y} = x^T A^T$$

**Proposition 3.4** For the special case in which the scalar α is given by the quadratic form $\alpha = x^T A x$ where $x$ is $n \times 1$, $A$ is $n \times n$, and A does not depend on $x$, then
$$\frac{\partial \alpha}{\partial x} = x^T(A + A^T).$$

### 4. MAIN RESULTS

In this section, we just integrate matrix that has dimension 2 by 2 and 3 by 3, Fabiounse number used to find the formula of integration. The definition of matrix is same as in the first section, in the following we see some proposition which determine matrix integration in special case. Proposition 1 and 2 deal with $3^{rd}$ dimension but 3 and 4 deal with $2^{nd}$ dimension matrix.

Now if we assume $y = \alpha(x)$ where $y$ is vector with n element and $x$ is m element vector so

$$\int y_k dx_i = \begin{bmatrix} \int a_{11} x_1 dx_1 & \cdots & \int a_{1n} x_n dxn \\ \vdots & \ddots & \vdots \\ \int a_{m1} x_1 & \cdots & \int a_{mn} x_n \end{bmatrix} \quad . \quad (1)$$

**Theorem 4.1** If $y_k = \varphi(x_i)$, where $y$ is depending in $x$, then integration will be as follow:

$$\int y_k dx_i = \begin{cases} x_i \sum_{j=1}^{3} \frac{1}{(n-j)!} a_{kj} x_j & , i = 1 \\ x_i \sum_{j=1}^{3} \left(\frac{n-j}{i}\right)! a_{kj} x_j & , i = 2 \\ x_i \sum_{j=1}^{3} \frac{1}{(j-1)!} a_{kj} x_j & , i = 3 \end{cases} \quad . \quad . \quad (2)$$

Where $= 1,2,3$ , $i = 1,2,3$.

**Proof**

Since $y_k = \sum_{j=1}^{3} a_{kj} x_j$ ,then the integration of $y_k$ respect to $x_i$ is
$$\int y_k dx_i = \sum_{j=1}^{3} \int a_{kj} x_j \quad . \quad (3)$$

Then if we integrate $y_k$ with respect to $x_1$ and return the equation from matrix to equation the integration formula will be as follow:
$$\int y_k dx_i = x_1 \sum_{j=1}^{3} \frac{1}{(n-j)!} a_{kj} x_j$$
and, if we integrate $y_k$ resp ect to $x_2$ can be written as:

$$\int y_k dx_i = x_2 \sum_{j=1}^{3} \left(\frac{n-j}{2}\right)! \, a_{kj} x_j$$

Also, if we integrate $y_k$ respect to $x_3$ can be written as:

$$\int y_k dx_i = x_3 \sum_{j=1}^{3} \frac{1}{(j-1)!} a_{kj} x_j$$

Hence by above integrations formulas we get(2).

**Corollary 4.1.1** Let $\beta$ be scalar that defend as follow $\beta = x^T A x$ then the integration is:

$$\int \beta dx_i = \begin{cases} x_i \sum_{k=1}^{3} x_k \sum_{j=1}^{3} \frac{1}{n-F(j-1)-[k-(k-1)]!} a_{kj} x_j & , i = 1 \\ x_i \sum_{k=1}^{3} x_k \sum_{j=1}^{3} \frac{1}{n-F(|j-i|)+(i-k)} a_{kj} x_j & , i = 2 \\ x_i \sum_{k=1}^{3} x_k \sum_{j=1}^{3} \frac{1}{F(j)+(n-j)^{n-k}} a_{kj} x_j & , i = 3 \end{cases} \quad . \quad (4)$$

where $A$ is 3 by 3 matrix and $i = 1,2,3$.

**Proof:**
Since $\beta = x_i \sum_{k=1}^{3} x_k \sum_{j=1}^{3} a_{kj} x_j$ , when i=1,2,3 then the integration will be
$\int \beta dx_i = x_i \int \sum_{k=1}^{3} x_k \sum_{j=1}^{3} a_{kj} x_j$

When we integrate $\beta$ with respect to $x_1$

$$\int \beta dx_i = x_1 \sum_{k=1}^{3} x_k \sum_{j=1}^{3} \frac{1}{n-F(j-1)-[k-(k-1)]!} a_{kj} x_j$$

and, if we integrate $\beta$ respect to $x_2$ it will be**:**

$$\int \beta dx_i = x_2 \sum_{k=1}^{3} x_k \sum_{j=1}^{3} \frac{1}{n-F(|j-i|)+(i-k)} a_{kj} x_j$$

also, if we integrate $\beta$ respect to $x_3$ it will be

$$\int \beta dx_i = x_3 \sum_{k=1}^{3} x_k \sum_{j=1}^{3} \frac{1}{F(j)+(n-j)^{n-k}} a_{kj} x_j$$

Thus from all integrations formulas, we find that (4).

**Theorem 4.2** If $y_k = \varphi(x_i)$, where $y$ is depending in $x$, then integration will b

$$\int y_k dx_{ik} = \begin{cases} x_i \sum_{j=1}^{2} \frac{j}{n} a_{kj} x_j & , i = 1 \\ x_i \sum_{j=1}^{2} \frac{1}{i} a_{kj} x_j & , i = 2 \end{cases} \quad . \quad . \quad . \quad (5)$$

where $k = 1,2$.

**Proof:**
Since $y_k = \sum_{j=1}^{2} a_{kj} x_j$ , then the integration of $y_k$ respect to $x_i$ Is
$$\int y_k dx_i = \sum_{j=1}^{2} \int a_{kj} x_j$$

If the integrate with respect to $x_1$:
$$\int y_k dx_{ik} = x_1 \sum_{j=1}^{2} \frac{j}{n} a_{kj} x_j \quad . \quad . \quad . \quad (6)$$

But if it is with respect to $x_2$ it is :
$$\int y_k dx_{ik} = x_2 \sum_{j=1}^{2} \frac{1}{i} a_{kj} x_j \quad . \quad . \quad . \quad (7)$$

So that from (6) and (7), we have (5).

**Corollary 4.2.1** Let $\beta$ be scalar that defend as follow $\beta = x^T A x$ then the integration is:

$$\int \beta dx_i = \begin{cases} x_i \sum_{k=1}^{2} x_k \sum_{j=1}^{2} \frac{1}{(n+1)-F(jk-i)} a_{kj}x_j & , i = 1 \\ x_i \sum_{k=1}^{2} x_k \sum_{j=1}^{2} \frac{1}{2F(k+j)-(2-j)k} a_{kj}x_j & , i = 2 \end{cases} \qquad . \qquad . \qquad (8)$$

where $A$ is matrix of order 2.

**Proof:** Since $\beta = x_i \sum_{k=1}^{2} x_k \sum_{j=1}^{2} a_{kj}x_j$ ,when $i = 1,2,3$ then the integration will be
$\int \beta dx_i = x_i \int \sum_{k=1}^{2} x_k \sum_{j=1}^{2} a_{kj}x_j$

If the integration was for $x_1$
$$\int \beta dx_i = x_i \sum_{k=1}^{2} x_k \sum_{j=1}^{2} \frac{1}{(n+1)-F(jk-i)} a_{kj}x_j \qquad . \qquad . \qquad . \qquad (9)$$

And when the integration was for $x_2$
$$\int \beta dx_i = x_i \sum_{k=1}^{2} x_k \sum_{j=1}^{2} \frac{1}{2F(k+j)-(2-j)k} a_{kj}x_j . \qquad . \qquad . \qquad . \qquad (10)$$

Therefore, from (9) and (10), we get (8).

## 5. APPLICATIONS

**Example 5.1:** Consider
$$y_1 = 2x_1 + 3x_2 + 4x_3$$
$$y_2 = 3x_1 + 5x_2 + 6x_3$$
$$y_3 = 7x_1 + 8x_2 + 6x_3$$

then, the system is $3 \times 3$ equation and y depending on the x so the matrix will be as follow after integration it with respect to $x_1$ and where $k = 1,2$ and 3:

$$\int y_k \, x_1 = \begin{bmatrix} x_1^2 & 3x_1x_2 & 4x_1x_3 \\ \frac{3}{2}x_1^2 & 5x_1x_2 & 6x_1x_3 \\ \frac{7}{2}x_1^2 & 8x_1x_2 & 6x_1x_3 \end{bmatrix}$$

But by **Theorem (4.1)** when $i = 1$ for $y_1$ with respect to $x_1$ we get

$$\int y_1 dx_1 = x_1 \sum_{j=1}^{3} \left( \frac{1}{(3-j)!} \right) a_{kj}x_j$$
$$= x_1 \left( \frac{1}{3-1} 2x_1 + \frac{1}{3-2} 3x_2 + \frac{1}{(3-3)!} 4x_3 \right)$$
$$= x_1^2 + 3x_1x_2 + 4x_1x_3.$$

**Example 5.2:** If
$$y_1 = 3x_1 + 4x_2$$
$$y_2 = 2x_1 + 3x_2,$$
then by corollary 2 we integral $\beta$ with respect to $x_1$ where it is $\beta = x^T A x$ so

$$\int \beta dx_1 = x_1 \sum_{k=1}^{2} x_k \sum_{j=1}^{2} \frac{1}{(n+1)-F(jk-i)} a_{jk}x_j$$

$$= x_1 \left[ \frac{3}{3-F(0)} x_1^2 + \frac{4}{3-F(3)} x_2^2 \right]$$
$$= x_1^3 + 4x_2^2 x_1$$

### 6. CONCLUSION

To sum Above, after using Fibonacci numbers to integrate square matrices of order 2 and 3 in the Sukhvinder algorithm, we noticed that the process of clarifying images might be faster.